\documentclass[12pt]{article}
\usepackage{amsmath,amssymb,amsthm}
\usepackage{fullpage}

\newtheorem{thm}{\textbf{Theorem}}

\newcommand{\dis}{\displaystyle}

\begin{document}

\title{Ramanujan's Inverse Elliptic Arc Approximation }
\author{Mark B. Villarino\\
Depto.\ de Matem\'atica, Universidad de Costa Rica,\\
2060 San Jos\'e, Costa Rica}
\date{\today}

\maketitle

 \begin{abstract}
 We suggest a continued fraction origin to \textsc{Ramanujan}'s approximation to $(\frac{a-b}{a+b})^{2}$ in terms of the arc length of an ellipse with semiaxes $a$ and $b$.\end{abstract}

Entry 14 of Chapter 38 (page 541) in \textsc{Berndt}'s edition of Volume 5 of \emph{Ramanujan's Notebooks} \cite{B} states (in part): \textit{\begin{quotation}\label{E}``Let the perimeter of ellipse $=\pi (a+b)(1+h),$ then\begin{equation}\fbox{$\dis \left(\frac{a-b}{a+b}\right)^{2}=4h-\frac{3h^{2}}{2+\sqrt{1-3h}}$}\end{equation}
\\
\indent very nearly . . ."\end{quotation}
}
Professor \textsc{Berndt}'s subsequent discussion of this entry merits comment on two points:\begin{enumerate}
  \item If one writes \begin{equation}
\label{lambda }
\lambda^{2}(h)\approx \dis 4h-\frac{3h^{2}}{2+\sqrt{1-3h}}
\end{equation}then \textsc{Berndt} discusses the (asymptotic) accuracy of the \textbf{\emph{inverse}} approximation\begin{equation}
h\approx h(\lambda^{2}),
\end{equation} i.e., how well the function of $\lambda^{2}$ approximates $h$, instead of the stated \textbf{\emph{direct}} approximation, i.e., how well the approximative function of $h$ in (1) approximates $\lambda^{2}.$
  \item \textsc{Berndt} provides \textbf{\emph{no}} suggestion as to how \textsc{Ramanujan} \textbf{\emph{might have proved (or discovered) the approximation }}$(1)$.
 \end{enumerate}

We therefore offer the following treatment of Entry 14 in which we deal with both points.  Our starting point is \textsc{Ivory}'s well-known series representation for the perimeter of an ellipse with semiaxes $a$ and $b$ (\cite{Per}, p. 146):\begin{thm}If \begin{equation}
\lambda:=\dis\frac{a-b}{a+b}, 
\end{equation}and if $L$ denotes the perimeter of the ellipse \begin{equation}
\dis\frac{x^{2}}{a^{2}}+\frac{y^2}{b^{2}}=1
\end{equation}then\begin{equation}
\label{Per }
L=\pi (a+b)\left\{1+\left(\frac{1}{2}\right)^{2}\lambda^{2}+\left(\frac{1\cdot 1}{2\cdot 4}\right)^{2}\lambda^{4}+\left(\frac{1\cdot 1\cdot 3}{2\cdot 4\cdot 6}\right)^{2}\lambda^{6}+\cdots\right\}.
\end{equation}
\end{thm}\hfill$\Box$

We deal with the first point.
\\

Writing\begin{equation}
L:=\pi(a+b)(1+h)
\end{equation}and using (6) to solve for $h$ we obtain the equation \begin{equation}h
\equiv h(\lambda)=\frac{1}{4}\lambda^{2}+\frac{1}{4^{3}}\lambda^{4}+\frac{1}{4^{4}}\lambda^{6}+\frac{25}{4^{7}}\lambda^{8}+\frac{49}{4^{8}}\lambda^{10}+\frac{441}{4^{10}}\lambda^{12}+\frac{1089}{4^{11}}\lambda^{14}+\frac{184041}{4^{15}}\lambda^{16}+\frac{511225}{4^{16}}\lambda^{18}+\cdots.
\end{equation}  If we \textbf{\emph{revert}} this series we obtain the \textbf{\emph{true}} expansion:\begin{equation}
\lambda^{2}\equiv\left(\frac{a-b}{a+b}\right)^{2}=4h-h^{2}-\frac{1}{2}h^{3}-\frac{5}{8}h^{4}-\frac{17}{16}h^{5}-\frac{\mathbf{273}}{128}h^{6}-\frac{609}{128}h^{7}-\frac{23391}{2048}h^{8}+\cdots
\end{equation}On the other hand, if we develop the \textbf{RHS} of \textsc{Ramanujan}s approximation (1) into powers of $h$ we obtain the \textbf{\emph{approximative}} expansion: \begin{equation}
\lambda^{2}\equiv\left(\frac{a-b}{a+b}\right)^{2}\approx4h-h^{2}-\frac{1}{2}h^{3}-\frac{5}{8}h^{4}-\frac{17}{16}h^{5}-\frac{\mathbf{269}}{128}h^{6}-\frac{1163}{256}h^{7}-\frac{10657}{1024}h^{8}-\cdots
\end{equation}Therefore we conclude that,$$(\text{\textbf{True} $\lambda^{2}) -$ (\textbf{Approximate}}\  \lambda^{2}) =-\frac{h^{6}}{32}-\frac{55}{256}h^{7}-\frac{2077}{2048}h^{8}-\cdots$$and thus, asymptotically, the following error estimate is valid:\begin{thm} \textsc{Ramanujan}'s approximation \textbf{\emph{overestimates}} the true value by about $\mathbf{\dfrac{1}{32}h^{6}}.$\end{thm}\hfill$\Box$ 

The accuracy is quite impressive.
\\

Now we deal with the second point.
\\

If we expand the \textbf{RHS} of the \textbf{\emph{true}} value $(9)$ into a \textbf{\emph{continued fraction}} we obtain\begin{equation}
\left(\frac{a-b}{a+b}\right)^{2}=4h-\cfrac{h^{2}}{ 1
                                   -\cfrac{\frac{1}{2}h}{ 1
                                   -\cfrac{\mathbf{\frac{3}{4}h}}{ 1
                                   -\cfrac{\mathbf{\frac{3}{4}h}}{ 1
                                   -\cfrac{\frac{3\frac{2}{9}}{2}h}{ 1-\cdots}}}}}
\end{equation}

We suggest that starting with the third numerator on, \textsc{Ramanujan} \textbf{\emph{took all of the numerators}} $\mathbf{=\dfrac{3}{4}h}$.  Then\begin{equation}
\left(\frac{a-b}{a+b}\right)^{2}\approx4h-\frac{h^{2}}{ 1
                                   -\cfrac{\frac{1}{2}h}{ 1
                                   -\cfrac{\mathbf{\mathbf{\frac{3}{4}h}}}{ 1
                                   -\cfrac{\mathbf{\mathbf{\frac{3}{4}h}}}{ 1
                                   -\cfrac{\mathbf{\frac{3}{4}h}}{\dis 1-\cdots}}}}}.
\end{equation}  Let us define $B$ by:\begin{equation}
\left(\frac{a-b}{a+b}\right)^{2}\approx4h-\frac{h^{2}}{\dis 1
                                   -\cfrac{\frac{1}{2}h}{\underbrace{\dis 1
                                   -\cfrac{\mathbf{\frac{3}{4}h}}{ 1
                                   -\cfrac{\mathbf{\frac{3}{4}h}}{ 1
                                   -\cfrac{\mathbf{\frac{3}{4}h}}{ 1-\cdots}}}}_{\equiv B}}}\end{equation}i.e., $B$ is the continued fraction\begin{equation}
B:= 1
                                   -\cfrac{\mathbf{\frac{3}{4}h}}{ 1
                                   -\cfrac{\mathbf{\frac{3}{4}h}}{ 1
                                   -\cfrac{\mathbf{\frac{3}{4}h}}{ 1-\cdots}}}.
\end{equation}To evaluate $B$ we observe that $B$ satisfies the equation\begin{align*}
&B=1-\frac{\frac{3}{4}h}{B}\\
&\Rightarrow B^{2}-B+\frac{3}{4}h=0\\
&\Rightarrow B=\frac{1+\sqrt{1-3h}}{2}
\end{align*}since $B=1$ when $h=0.$ Therefore, (13) becomes: \begin{align*}\left(\frac{a-b}{a+b}\right)^{2}&\approx4h-\frac{h^{2}}{\dis 1
                                   -\cfrac{\frac{1}{2}h}{\dfrac{1+\sqrt{1-3h}}{2}
                                   }}\\
                                   &=4h-\frac{3h^{2}}{2+\sqrt{1-3h}}                                   
\end{align*}and we have recovered (or \textbf{\emph{discovered!}}) \textsc{Ramanujan}'s approximation (1).

This derivation also explains why \textsc{Ramanujan}'s approximation is so accurate, since the continued fraction of the true value coincides with that of the approximative value for the first \textbf{\emph{four}} convergents so that one knows that the error is that of the fifth convergent.

\textsc{Ramanujan}'s mastery of continued fractions is justly famous, and we respectfully suggest that he discovered the approximation (1) by a method that follows the march of our presentation.

\subsubsection*{Acknowledgment}
Support from the Vicerrector\'{\i}a de Investigaci\'on of the 
University of Costa Rica is acknowledged.

\end{document}